\documentclass[12pt]{amsart}

\usepackage{amsmath,amssymb}

\newtheorem{thm}{Theorem}%[section]
\newtheorem{prop}[thm]{Proposition}%[section]
\newtheorem*{rem}{Remark}

\newcommand{\Z}{{\mathbb Z}}

\newcommand{\fp}{{\mathbb F}_p}
\newcommand{\fpn}{{\mathbb F}_{p^n}}
\newcommand{\fq}{{\mathbb F}_q}
\newcommand{\fqn}{{\mathbb F}_{q^n}}
\newcommand{\fqnp}{{\mathbb F}_{q^{n'}}}

\newcommand{\pK}{{\mathfrak p}}
\newcommand{\pE}{{\mathfrak Q}}
\newcommand{\pKa}{{\mathfrak p}_a}
\newcommand{\pKp}{{\mathfrak p}'}
\newcommand{\pLp}{{\mathfrak P}}

\newcommand{\OLp}{{{ O}_{L'}}}
\newcommand{\OKp}{{{ O}_{K'}}}
\newcommand{\OK}{{{ O}_{K}}}
\newcommand{\OOE}{{{O}_{E}}}

\newcommand{\gal}{\operatorname{Gal}}
\newcommand{\unramified}{\mathbb{P}'}
\newcommand{\frob}{\operatorname{Frob}_p}
\newcommand{\frobq}{\operatorname{Frob}_q}
\newcommand{\ord}{\operatorname{ord}}

\newcommand{\artin}[2]{\left( \frac{#1}{#2} \right)}

\newcommand{\fixme}[1]{}

\newcommand{\fixmelater}[1]{}

\newcommand{\C}{{\mathbb C}} 

\title[The Chebotarev density theorem --- incomplete intervals]{The Chebotarev density theorem for function fields --- incomplete intervals}

\author{P\"ar Kurlberg, Lior Rosenzweig}
\urladdr{www.math.kth.se/\~{ }kurlberg}
\address{Department of Mathematics, KTH Royal Institute of Technology,
SE-100 44 Stockholm, Sweden}
\email{kurlberg@math.kth.se}

\address{Unit of Mathematics, Afeka Tel Aviv College of Engineering, Mivtza Kadesh 38, Tel Aviv, Israel}
\email{liorr@afeka.ac.il}

\thanks{P.K. was partially supported by
the
Swedish Research Council (2016-03701).}

\date{July 3, 2020}
\begin{document}
\begin{abstract}
  We prove a P\'olya-Vinogradov type variation of the the Chebotarev
  density theorem for function fields over finite fields valid for
  ``incomplete intervals'' $I \subset \fp$, provided
  $(p^{1/2}\log p)/|I| = o(1)$.  Applications include density results
  for irreducible trinomials in $\fp[x]$, i.e. the number of
  irreducible polynomials in the set
  $\{ f(x) = x^{d} + a_{1} x + a_{0} \in \fp[x] \}_{a_{0} \in I_{0},
    a_{1}\in I_{1}}$ is $\sim |I_{0}|\cdot |I_{1}|/d$ provided
  $|I_{0}| > p^{1/2+\epsilon}$, $|I_{1}| > p^{\epsilon}$, or
  $|I_{1}| > p^{1/2+\epsilon}$, $|I_{0}| > p^{\epsilon}$, and
  similarly when $x^{d}$ is replaced by any monic degree $d$
  polynomial in $\fp[x]$.  Under the above assumptions we can also
  determine the distribution of factorization types, and find it to
  be consistent with the distribution of cycle types of 
  permutations in the symmetric group $S_{d}$.
\end{abstract}
%\today
%\timestamp
\maketitle

\section{Introduction}
\label{sec:introduction}

The distribution of primes, and more generally the distribution of
factorization types, in ``short intervals'' in the setting of function
fields over finite fields has received considerable attention
\cite{cohen-uniform,cohen-odoni77,BBR,BB,kr-shortintervals}.
For example, in \cite{BBR}, prime equidistribution for the family
$\{ f(x) + bx + a\}_{a,b\in \fp}$ was shown for
$f \in \fp[x]$ any monic degree $d$ polynomial (for $p$ large.)
For ``very short'' intervals, i.e., one parameter families of the form
$\{f(x) +a \}_{a \in \fp}$, prime equidistribution does not hold for
all $f$.  However, for $f$ suitably ``generic'', prime
equidistribution does in fact hold for very short intervals. In
\cite{kr-shortintervals} it was shown that given a monic degree $d$
``Morse polynomial'' $f(x) \in \fp[x]$ (i.e., that
$|\{ f(\xi) : f'(\xi) = 0 \} | = d-1$; note this holds for generic
polynomials, cf. Section~\ref{sec:proof-theorem})
\begin{equation}
  \label{eq:irreducible-count}
\left|\left\{ a \in \fp : f(x) + a \text{ irreducible }\right\}\right|
  =
p/d + O_d(p^{1/2}).
\end{equation}
More generally, the distribution of  factorization types of
$f(x)+a$ can also be determined.  Writing
$f(x) = \prod_{i=1}^{l} f_{i}(x)$ with $f_{i} \in \fp[x]$ all
irreducible and letting $d_{i} = \deg(f_{i})$ we may after rearranging
assume that $1 \leq d_{1} \leq d_{2} \ldots \leq d_{l}$; the
{\em factorization type} (or decomposition type) of $f$ is then given by
$(d_{1}, \ldots, d_{l})$.
The distribution of factorization types of $f(x)+a$, for $f$ Morse,
is consistent (up to an error of size $O_{d}(p^{-1/2}))$ with the
distribution of cycle types of permutations in $S_{d}$, the
symmetric group on $d$ letters, with respect to the Haar measure.
(E.g., for $\sigma = (12) \in S_{3}$, write out all trivial cycles,
i.e., $\sigma = (12)(3)$; the {\em cycle type} of $\sigma$ is then
$(1,2)$ if we order according to cycle lengths.)

The connection between factorization types and group theory proceeds
via Galois theory and the function field version of the Chebotarev
density theorem, made effective by Weil's proof of the Riemann
hypothesis for curves.  The key point is that for $f(x)$ Morse,
$\gal(f(x)+t / \fp(t)) \simeq S_{d}$ and the factorization type for
$f(x)+a$ can be read off from the cycle type of the Frobenius
class
%, regarded as a subset of $S_{d}$, 
at the prime ideal $(t-a) \subset \fp[t]$.  In particular, $f(x)+a$
being irreducible is equivalent to the Frobenius conjugacy class at
the prime
$(t-a)$ being generated by a $d$-cycle, and the proportion of
$d$-cycles in $S_{d}$ is $1/d$, hence the density $1/d$ in
(\ref{eq:irreducible-count}).

The purpose of this paper is to show that equidistribution of
factorization types also holds for significantly smaller subsets,
namely for ``incomplete intervals'' $I \subset \fp$, as long as 
$(p^{1/2}\log p) / |I|$ is small.
In fact, in spirit of the P\'olya-Vinogradov inequality, we will develop
a version of the Chebotarev density theorem for incomplete intervals,
allowing us to determine the distribution of Artin symbols, and
thus resolve finer invariants than factorization types when
$\gal(f(x)+t /  \fp(t))$ is not the full symmetric group.  (Note that the
cycle type is in general not enough to determine the
conjugacy class,  e.g., the three-cycles $(123)$ and $(132)$ are not
conjugate in the alternating group $A_{3}$.)

Before stating our main result we introduce some notations.
Let $p$ be a (large) prime, let $K :=
\fp(t)$, let $L/K$ be a finite normal and separable 
extension with Galois group $\gal(L/K)$, and let $\OK = \fp[t]$ denote
the ring of integers in $K$.
%, and let $\OL$ denote the integral
% closure of $\OK$ in $L$.
Given a prime ideal $\pK \subset \OK$ that does not ramify in $L$, let
$\artin{L/K}{\pK} $ denote the Artin symbol, a certain conjugacy class
in $\gal(L/K)$.  (For further details and definitions,
cf. \cite[Ch.~6]{FJ}.)  It will be convenient to use the
convention that any prime $\pK \subset \OK$ appearing in an Artin symbol
$\artin{L/K}{\pK}$ is implicitly assumed to be unramified.

By an {\em incomplete interval} $I$ in $\fp$ we mean a
set of the form $I = [M,M+1, \ldots, M+N] \subset \fp$ for
$M \in \fp$ and $N \in \mathbb{N}$ (in fact, our method applies to arithmetic
progressions of the form $\{ A \cdot i + B \}_{i \in I}$ where
$A,B\in \fp$ and $I$ is an incomplete interval.)
Let $\frob$ denote the Frobenius substitution $\alpha \to \alpha^{p}$.
Define $n \in \Z^+$ such that $\fpn = L \cap \overline{\fp}$ is the
field of constants in $L$, and let $m = [L: K \fpn]$; note that in the
``geometric case'', i.e., when $n=1$, we have $m= |\gal(L/K)|$.
Further, given
$a \in \fp$, let $\pKa \subset \fp[t]$ denote the prime ideal
generated by $(t-a)$.
%
% Given an incomplete interval $I \subset
%  \fp$ let
% $$
% I' := \{ a \in I : \pKa \text{ does not ramify in $L$ } \};
% $$
% note that $|I'| = |I| + O_{[L:K]}(1)$.

\begin{thm}
\label{thm:incomplete-chebotarev}  
Let $C \subset \gal(L/K)$ be a conjugacy class.
If 
$\tau|_{\fpn} = \frob|_{\fpn}$ for all $\tau \in C$, then
  $$
  \left|\left\{ a \in I : \artin{L/K}{\pKa} = C \right\}\right| =
  \frac{|C|}{m} \cdot |I| + O_{[L:K]}({p^{1/2}} \log p),
  $$
  On the other hand,
  if 
$\tau|_{\fpn} \neq \frob|_{\fpn}$ for all $\tau \in C$, then we have
$$ \left| \left\{ a \in I : \artin{L/K}{\pKa} = C \right\}\right|=0.$$
\end{thm}
The result only gives non-trivial information for $|I|$ slightly
larger than $p^{1/2} \log p$.  
Using the P\'olya-Vinogradov method of ``completing the sum'', the
result follows easily  from our key technical
result, 
Proposition~\ref{prop:twisted-chebotarev},
namely square root cancellation for certain complete sums
twisted by additive characters.
We remark that by using smoothing
(cf. \cite[Theorem~1.1]{fouvry-etal-short-trace-sums}) the error in
Theorem~\ref{thm:incomplete-chebotarev} is easily improved to
$O_{[L:K]}({p^{1/2}} \log (10^{5}|I|/p^{1/2}))$, hence giving
asymptotics as long as $p^{1/2}/|I|$ tends to zero.  We leave the
details to the interested reader.

\begin{prop}
  \label{prop:twisted-chebotarev}
  Let $K = \fp(t)$ and let $L/K$ be a normal and separable extension with Galois
  group $\gal(L/K)$, and let $C$ denote a conjugacy class in
  $\gal(L/K)$.  If $\psi : \fp \to \C^{\times}$ is a non-trivial
  additive character then
$$
  \left|
    \sum_{ a \in \fp : \artin{L/K}{\pKa} = C }{}\hspace{-.78cm}^{'} \psi(a) 
    \right|
\ll_{[L:K]} p^{1/2},
$$
with $\sum'$ denoting the sum restricted to $a \in \fp$ such that
$\pKa$ is in addition unramified in $L$.
\end{prop}
We remark that the sum is empty unless $\tau|_{\fpn} = \frob|_{\fpn}$
for all $\tau \in C$.

\subsection{Applications}
\label{sec:applications}

We now give some explicit examples of families of polynomials for
which we can determine the corresponding Galois group. This together
with Theorem~\ref{thm:incomplete-chebotarev} shows the existence of
primes in incomplete intervals.  The method proceeds by finding the
distribution of Artin symbols, hence also determines the distribution
of factorization types, but for simplicity we state it only for the
density of irreducible polynomials.

\begin{thm}
  \label{thm:primes-in-incomplete-interval}
Let $p > d+1$ be a (large) prime, and let $f(x) = x^{d}  + \sum_{i=0}^{d-1}
a_{i} x^{i} \in \fp[x]$ be a monic 
polynomial of degree $d$.
There
exists subsets $B_{1}, B_{2}, B_{3} \subset \fp$ such that $|B_{i}| =
O_{d}(1)$ for $i=1,2,3$ with the following properties.  Given any
incomplete interval $I \subset \fp$, we have:
\begin{enumerate}
\item For $a_{1} \in \fp \setminus B_{1}$,
  $$
\left|\left\{ \alpha \in I : f(x) + \alpha  \text{ irreducible } \right\} \right| =
|I|/d + O_{d}(p^{1/2} \log p).
$$
\item
For $a_{0} \in \fp \setminus B_{2}$,
$$
|\{ \alpha \in I : f(x) + \alpha \cdot x \text{ irreducible } \} | =
|I|/d + O_{d}(p^{1/2} \log p).
$$
  
\item For $a_{1} \neq 0$, $a_{0} \in \fp \setminus B_{3}$ and any integer
  $m \in [2,d-1]$, 
  $$
|\{\alpha \in I : f(x) + \alpha \cdot x^{m}    \text{ irreducible } \} | =
|I|/d + O_{d}(p^{1/2} \log p).
$$
\end{enumerate}

\end{thm}

An immediate application is a
deterministic time algorithm for 
constructing irreducible degree $d$ polynomials $g(x) \in \fp[x]$ of
quite general shapes, in particular of very sparse form.  (The
complexity in terms of $d$ is quite bad so we will not make it
explicit.)  Given any (say monic) degree $d$ polynomial $f$, consider
the family $f(x) + bx +a$.  Trying at most $O_{d}(1)$ values of $b$,
and $O_{d}(p^{1/2} \log p)$ values of $a$ yields an irreducible
polynomial; each such irreducibility test can be done in
$O_{d}( \log p)$ arithmetic operations in $\fp$ (say using Rabin's
test).  Hence we can produce an irreducible $g(x)$ in
$O_{d}(p^{1/2} \log^{2} p)$ $\fp$-operations, or even
$O_{d}(p^{1/2} \log p)$ operations using the earlier described
smoothing improvement; note that for existence for irreducibles it is
enough to take $|I| \gg_{d} p^{1/2}$.  This is  to be 
compared with Shoup's algorithm
\cite{shoup-deterministic-irreducibility} which requires
$O_{d}(p^{1/2} \log^{3} p)$ $\fp$-operations.  Ignoring polynomial
factors in $\log \log p$, and using fast (FFT-based) $\fp$-arithmetic,
the bit operation complexities can be obtained by multiplying the
above bounds by $\log p$.

\begin{rem}
  Given polynomials $f,g \in \fp[x]$ with $\deg(g) < \deg(f) = d$ we can under
  fairly weak assumptions on $f,g$, namely that the ratio $f/g$ is
  Morse, show that (cf. Section~\ref{sec:crit-morse-polyn})
  $$
|\{\alpha \in I : f(x) + \alpha \cdot g(x)  \text{ irreducible } \} | =
|I|/d + O_{d}(p^{1/2} \log p).
$$
\end{rem}

For $f$ a degree $d$ Morse polynomial, and distinct
$h_{1}, \ldots, h_{k} \in \fp$, the Galois group of the compositums of
the fields generated, over $\fp(t)$, by the polynomials
$ f(x) + h_{1} + t, \ldots, f(x) + h_{k} +t$ is maximal, i.e.,
isomorphic to $S_{d}^{k}$ (cf. \cite{kr-shortintervals}).
Theorem~\ref{thm:incomplete-chebotarev} then gives equidistribution of
cycle types of
$\{(f(x) + h_{1} + a, \ldots, f(x) + h_{k} +a)\}_{a \in I}$ inside
$S_{d}^{k}$, provided $(p^{1/2}\log p)/|I| = o(1)$.  Thus, for $f$
Morse we immediately obtain cancellation for function field analogs of
Moebius and Chowla type sums; with $\mu$ denoting the function field
Moebius $\mu$ function, we have
$$
\sum_{a \in I} \mu( f + h_{1} + a) \cdot \ldots \cdot \mu( f + h_{k} +
a)
\ll_{k,d} p^{1/2} \log p
$$
(cf. \cite{kr-shortintervals,carmon-rudnick-moebius-chowla,KR} for
results valid for various longer intervals.)

Another application is asymptotics for function field anologs, for
incomplete intervals ``centered'' at $f$, assuming
$(p^{1/2} \log p)/|I| = o(1)$ and $f$ Morse, of shifted divisors
sums (e.g.  $\sum_{a \in I} d_{r}(f+a) d_{r}(f+1+a)$ where $d_{r}$ is
the $r$-th divisor function) and the Titchmarsh divisor problem (e.g.
$\sum_{a \in I} 1_{\operatorname{Prime}}(f+a) d_{r}(f+1+a)$);
asymptotics for these sums over ``long'' intervals, while allowing for
very general shifts, were determined in \cite{ABR}.

We next give an example where the cycle type distribution is not
enough to determine the distribution of Artin symbols.
Let $p \equiv 1 \mod 3$ be a (large) prime, and let
$h_{1}, \ldots, h_{k} \in \fp$ be distinct.  Define
$f_{i}(x) = x^{3} + h_{i} + t$ and let $L_{i}$ denote the splitting
field of $f_{i}$ over $\fp(t)$, and let $L^{k}$ denote the compositum
of the fields $L_{1},\ldots, L_{k}$.  Then
(cf. \cite[Proposition~8]{K09}) 
$\gal(L_{i}/\fp(t)) \simeq A_{3}$ for $i=1,\ldots,k$, and
$\gal(L^{k}/\fp(t)) \simeq A_{3}^{k}$.  With $\pK_{a} = (t-a)$,
Theorem~\ref{thm:incomplete-chebotarev} implies that
the Artin symbols $\{ \artin{L^{k}/\fp(t)}{\pK_{a}} : a \in I \}$
equidistribute in $A_{3}^{k}$, with relative error of size $(p^{1/2}
\log p) /|I|$.

\subsection{Discussion}

For the full interval $I=\fp$ and $n=1$ (the geometric case)
Chebotarev's density theorem for function fields, with error term of
size $O_{[L:K]}(p^{1/2})$, easily follows from the work of Reichardt
\cite{reichardt-chebotarev} together with Weil's celebrated proof of
the Riemann hypothesis for curves over finite fields
\cite{weil-artin-rh-book}.  (Reichardt's error term, for primes of
degree $v$, is of size $O(p^{\Theta \cdot v})$ as $v \to \infty$, where
$\Theta<1$ is the maximum of the real parts of the roots of the zeta function
of $L$.)
%
%
% goes back to Reichardt
% \cite{reichardt-chebotarev}.  However, as $\delta$
% might depend on $p$,  equidistribution for
% $\artin{L/\fp(t)}{\pK}$ only holds for  $\pK$ ranging over large degree
% primes in $\fp[t]$.
% Obtaining an error of size $O(p^{1/2})$
% is straightforward given Weil's celebrated proof of the Riemann
% hypothesis for curves over finite fields \cite{weil-artin-rh-book}.
%
%
The case $I= \fp$ and $n>1$ is due to Cohen and Odoni
%The case $n>1$ is more delicate and due to Cohen and Odoni
\cite{cohen-odoni77} (also cf.
\cite{jarden-chebotarev}).
As for incomplete intervals we are only aware of recent results by
Entin \cite{entin-incomplete} who determined the distribution of
factorization types for various families of polynomials whose
coefficients are allowed to vary over incomplete intervals (see
below).  
% In fact, though not stated, his methods can also determine
% the distribution of Artin symbols.
%worth noting
%is that their arguments only uses the Riemann hypothesis for curves.

%As for applications to irreducible polynomials, we note that
In \cite{shparlinski-on-primitive-elements-90} Shparlinski studied the
proportion of irreducible monic polynomials in $\fp[x]$, with
coeffecients in constrained to lie on points with integer coordinates
inside parallelepipeds, and showed that
\begin{multline*}
\left|
\left\{ (a_{0}, \ldots, a_{d-1}) \in I_{0} \times \cdots \times I_{d-1} : x^{d} +
\sum_{i=0}^{d-1} a_{i} x^{i} \text{ irreducible } \right\}
\right|
\\=
\frac{ \prod_{i=0}^{d-1} |I_{i}|}{d} + O( p^{d-1} \log^{d-1} p),
\end{multline*}
giving non-trivial information when
$\prod_{i=0}^{d-1} |I_{i}| > p^{d-1+\epsilon}$.  (More generally, he
also determined the distribution of factorization types.)
Further, in \cite{shparlinski-irreducible-trinomials}, he considered
sparser families, namely the set of irreducible trinomials in
$\fp[x]$, and showed that
$$
\bigg|\{ (a_{0},a_{1}) \in I_{0} \times I_{1} :
x^{d} + a_{1} x + a_{0}  \text{ irreducible } \} \bigg|
\sim|I_{0}|\cdot |I_{1}| /d
$$
provided $|I_{0}|,|I_{1}| > p^{1/4+\epsilon}$ and
$|I_{0}|\cdot |I_{1}| > p^{1+\epsilon}$, and similarly for trinomials
with any prescribed factorization type.
Theorem~\ref{thm:primes-in-incomplete-interval} easily implies
the same asymptotics, under the weaker conditions
$|I_{0}| > p^{\epsilon}$, $|I_{1}| > p^{1/2+\epsilon}$, or
$|I_{1}| > p^{\epsilon}$, $|I_{0}| > p^{1/2+\epsilon}$.

Entin \cite{entin-incomplete} determined the distribution of
factorization patterns for families of polynomials whose coefficients
vary over quite general sets $S \subset \fp^{d}$, with relative error
of size $\operatorname{irreg}(S)/p^{1/2}$ where
$\operatorname{irreg}(S)$ is related to the decay of Fourier
coefficients of the characteristic function of $S$.  E.g., for
$S = I_{0} \times \ldots \times I_{d-1}$,
$\operatorname{irreg}(S) \ll p^{d}\log^{d} p / (\prod_{i=0}^{d-1}
|I_{i}|)$ with relative error small if
$\prod_{i=0}^{d-1} |I_{i}| > p^{d-1/2+\epsilon}$.
%, a stronger assumption than Shparlinski's.  
His method also applies for sparser
families: for trinomials the assumption needed for small relative
error is $|I_{0}| \cdot |I_{1} | > p^{3/2 + \epsilon}$; for $f$ Morse
and the family $f(x) + a_{0}$ the relative error is small if
$(p^{1/2} \log p)/|I_{0}| = o(1)$, similar to the conditions in
Theorem~\ref{thm:primes-in-incomplete-interval}.

\subsection{Acknowledgements}
\label{sec:acknowledgements}
We thank A. Granville, E. Kowalski, Z. Rudnick and P. Salberger for helpful
discussions.  We also thank I. Shparlinski and O. Gorodetsky for valuable
comments on an early version of the paper,
in particular for suggesting using smoothing to improve the error term
in Theorem~\ref{thm:incomplete-chebotarev}.  We further
thank I. Shparlinski for pointing out the application of finding
irreducible polynomials in sparse families, thus leading to a
deterministic way to find irreducibles that is faster than Shoup's
algorithm.

\section{Deducing Theorem~\ref{thm:incomplete-chebotarev} from
  Proposition \ref{prop:twisted-chebotarev}}

To simplify the notation, let $F(a)=1$ if $\artin{L/K}{\pK} = C$,
and zero otherwise.  For $b \in \fp$, let $\psi_{b}(x)$ denote the
additive character $\psi_{b}(x) := e^{2\pi i bx/p}$.  With $1_{I}$
denoting the characteristic function of $I$, we ``complete the sum''
and write
$$
|\{ a \in I : \artin{L/K}{\pKa} = C \}|
=
\sum_{a \in I} F(a)
=
\sum_{a \in \fp} F(a) 1_{I}(a)
$$
$$
=
\sum_{a \in \fp} F(a)
\sum_{b \in \fp} 
\widehat{1_{I}}(b) \psi_{b}(a)
=
\sum_{b \in \fp} 
\widehat{1_{I}}(b)
\sum_{a \in \fp} F(a) \psi_{b}(a)
$$
where $\widehat{1_{I}}(b) := 1/p \sum_{c \in \fp} 1_{I}(c)
\psi_{b}(-c)$ is the Fourier transform of $1_{I}$.

Our main term will come from terms with $b=0$; by
\cite[Proposition~6.4.8]{FJ} (note that the genus of $L$ can be
bounded in terms of $[L:K]$) we find that 
\begin{multline}
    \label{eq:main-term1}
\frac{|I|}{p}
\sum_{a \in \fp} F(a)
=
\frac{|I|}{p} \cdot
\left(
  \frac{|C|}{m} p + O_{[L:K]} (p^{1/2})
\right)
\\=
 \frac{|C|}{m} |I| + O_{[L:K]} (|I|p^{-1/2})
\end{multline}

As for the error terms, taking integer representatives $b$ of the
elements in $\fp^{\times}$, such that $0 < |b| < p/2$, we find that
$\widehat{1_{I}}(b) \ll 1/b$.  Now, by
Proposition~\ref{prop:twisted-chebotarev} (note that for $b \not
\equiv 0 \mod p$, 
$\psi_{b}$ is non-trivial),  we have
\begin{equation}
  \label{eq:main-term2}
\sum_{b \in \fp^\times}
\left|
\sum_{a \in \fp} F(a) \psi_{b}(a)
\right|
\ll_{[L:K]} p^{1/2} \sum_{b=1}^{p/2} 1/b
\ll_{[L:K]} p^{1/2} \log p.
\end{equation}
As $|I| \leq p$, the error term in (\ref{eq:main-term2}) dominates the
error term in (\ref{eq:main-term1}), hence
$$
\frac{|I|}{p}
\sum_{a \in \fp} F(a)
=
\frac{|C|}{m} |I| + O_{[L:K]}(p^{1/2} \log p).
$$

\section{Additive characters and Artin-Schreier extensions}

We briefly recall how a non-trivial additive character $\psi$ on $\fp$
can be realized as a character $\psi_{K}$ on some cyclic galois group
$\gal(E/K)$, where  $K=\fp(t)$,  such that
%the Artin symbol should then satisfy
$$
\psi_K \left(
  \artin{E/K}{\pK_{a}}
\right)
= \psi(a),
\quad a \in \fp;
$$
and $\pK_{a}= (t-a) \subset \fp[t]$ denotes a degree one prime ideal.
(Note that $\deg(\pK)$, the degree of the prime ideal $\pK$ is defined
in terms of the cardinality of the residue field, namely
$|\OK/\pK| = p^{\deg(\pK)}$.  Also, since $E/K$ is cyclic we can
regard the Artin symbol $ \artin{E/K}{\pK_{a}}$ as an element in
$\gal(E/K)$ rather than a conjugacy class.)

The extension $E$ can be constructed as follows:
with $\xi$ denoting a root of the polynomial $f(x) = x^{p} - x -t$,
let $E=K(\xi)$.  Then $E/K$ is an Artin-Schreier extension of degree
$p$, such that $\gal(E/K) \simeq \fp^{+}$. The extension is unramified
except at infinity, where it is wildly ramified.

Specializing at $\pK_{a}=(t-a)$ for $a=0$, we find that $x^{p}-x$ splits
completely, and thus $\artin{E/K}{\pK_{a}}$ is the identity
element in $\gal(E/K)$.
Similarly, specializing at $\pK_{a}=(t-a)$ for $a \in \fp^\times$, we find that
$x^{p} - x - a$ is irreducible (over $\fp$), and thus there is exactly
one (unramified) prime $\pE_{a} \subset \OOE$ above $\pK_{a}$, and
$\OOE/\pE_{a} \simeq \mathbb{F}_{p^{p}}$.

We next turn to identifying the  Galois action of the
Artin symbol.
If $\xi$ is some fixed root of $x^{p}-x-t$, any Galois element
$\sigma \in \gal(E/K)$ acts via $\xi \to \xi + \alpha$ for some
$\alpha \in \fp$.  Moreover, given such a $\sigma$, we can recover
$\alpha$ from $\sigma(\xi) - \xi$.  Thus, if $\sigma$ is the
Artin map at $\pK_{a}=(t-a)$, we find that
$$
\alpha = \sigma(\xi) - \xi = \xi^{p} - \xi = a
$$

In conclusion, for $a \in \fp$ the Artin symbol $\artin{E/K}{\pK_{a}}$
acts as $\xi \to \xi + a$, hence there is a natural identification of
$\gal(E/K)$ with $\fp^{+}$ so that the image of $\pK_{a}$ is just
$a \in \fp^{+}$.  In particular, given a non-trivial additive character
of $\fp$ we can define a character $\psi_{K}$ on $\gal(E/K)$ such that
$$
{\psi}_{K} \left(
\left(
  \frac{E/K}{\pK_{a}}
\right)
\right)
=
\psi(a), \quad a \in \fp.
$$

\section{Proof of Proposition~\ref{prop:twisted-chebotarev}}
\label{sec:proof-proposition}
The main bulk of the argument is similar to the one used in
\cite[Ch~6.4]{FJ}; below we briefly
summarize the argument and only give details when we need to go beyond
their results. For easier comparison, we follow their notation and
write $q$ rather than $p$.

Let $K = \fq(t)$ and let $L/K$ be a normal extension with Galois group
$\gal(L/K)$.  
Since the constant in the error term is allowed to depend on $[L:K]$,
we may assume 
that $q=p$ is sufficiently large (say $q> [L:K])$ so that $L/K$ is
separable, as well as tamely ramified and thus
by the Hurwitz genus formula, that the genus  bound $g(L) \ll_{[L:K]} 1$ 
holds.
Let $C \subset \gal(L/K)$ be a conjugacy class.  We may
also assume that $\tau|_{\fqn} = \frobq|_{\fqn}$ for all $\tau \in C$,
otherwise the Artin symbol never meets $C$ and the sum is
empty.
%As we will mainly be concerned with unramifide degree one
%primes,
We will need some further notations: given field extensions
$L/K'/K$, let
\begin{multline*}
\mathbb{P}'_{k}(K')
:=
\{ \pK \subset \OKp : \text{ $\pK$ is a prime ideal,}
  \\ \text{ unramified over  $K$ and in $L$, }
\deg(\pK) = k \}.
\end{multline*}

Moreover, for extensions $L'/L/K'/K$, with $L'/L$ unramified, and a
conjugacy class $C' \subset \gal(L'/K')$, let
$$
C_{k}(L'/K', C')
:=
\left\{
\pKp \in \mathbb{P}'_{k}(K') : 
\artin{L'/K'}{\pKp} = C'
\right\}.
$$

Let $\fqn = L \cap \overline{\fq}$ denote the field of
constants in $L$, and let $m = [L:K \fqn]$.  
Choose some $\tau \in C$ and define $n' = n \cdot \ord(\tau)$, where
$\ord(\tau)$ denotes the order $\tau$ in $\gal(L/K)$; note that $n'$
does not depend on the choice of $\tau$ since all elements in $C$ have
the same order.
With $L' := L \fqnp$ we have $[L':K\fqnp] = [L:K\fqn] =
m$; also note that $L'/L$ is unramified.
By \cite[Lemma~6.4.4]{FJ}
there exists $\tau' \in \gal(L'/K)$ such 
that $\tau'|_{L} = \tau$ and $\tau'|_{\fqnp} = \frobq|_{\fqnp}$; 
moreover  $\ord(\tau') = \ord(\tau) \cdot n = n'$.

With $C' \subset \gal(L'/K)$ denoting the conjugacy class of $\tau'$,
we have (cf. \cite[Lemma~6.4.4]{FJ})
$$
C_{1}(L'/K, C') = C_{1}(L/K,C), 
$$
and thus
\begin{equation}
  \label{eq:CandCprimesums}
\sum'_{ a \in \fp : \artin{L/K}{\pKa} = C } \psi(a)   =
\sum_{\pK \in C_{1}(L/K, C)} \psi_{K}(\pK)
=
\sum_{\pK \in C_{1}(L'/K, C')} \psi_{K}(\pK)
\end{equation}

Let $K'_{\tau'} \subset L'$ denote the fixed field of $\tau'$.  Note
that $L'$ does not depend on $\tau'$, whereas $K'_{\tau'}$ does.  We will
need to keep track of this dependence and will therefore deviate slightly
from the notation in Fried-Jarden (in their argument it suffices to work with a
fixed $\tau'$, consequently they denote $K'$ for the fixed field of $\tau'$.)

From the proof of
\cite[Proposition~6.4.8]{FJ} it follows
that $K'_{\tau'} \cap \fqnp = K'_{\tau'} \cap \overline{\fq} = \fq$ ,
and $K'_{\tau'} \fqnp = L'$, as well as
$[K'_{\tau'}:K] = [L':K \fqnp] = [L:K\fqn] = m$.  It also follows
(in particular see how their Corollary~6.4.3 was used in the proof of
Proposition~6.4.8) 
$$
|C_{1}(L'/K,C')| 
= 
\frac{|C|}{[K'_{\tau'}:K]} \cdot |C_{1}(L'/K'_{\tau'}, \{ \tau'\})| 
= 
\frac{|C|}{m} \cdot |C_{1}(L'/K'_{\tau'}, \{ \tau'\})| 
$$
This equality suggests the existence of some $m$ to $1$ map; we will
show that this is indeed the case and then obtain a sum over certain
degree one primes in $\OKp$.

We first show that we can ``lift'' $\psi_{K}$ to a suitably invariant
character $\psi_{K'_{\tau'}}$, i.e., independent of the choice of
prime $\pKp$ above any prime $\pK$ occuring in the sums we wish to
estimate.
Recall that we have realised the additive character
$\psi : \fq \to \C^\times$ as a character $\psi_K$ on $\gal(E/K)$, for
$E/K$ an Artin-Schreier extension $E/K$, with the property that
$$
\psi_K \left( \artin{E/K}{\pK_a} \right) = \psi(a).
$$

As $[E:K] = q$ we have, for $q$ sufficiently large (which we may
assume to hold since the constant in error term is allowed to depend
on $[L:K]$), 
$\gal(K'_{\tau'}E/K) \simeq \gal(K'_{\tau'}/K) \times \gal(E/K)$, as
well as $\gal(K'_{\tau'}E/K'_{\tau'}) \simeq \gal(E/K)$, where the
latter isomorphism can be defined by restricting
$\gal(K'_{\tau'}E/K'_{\tau'})$ to act on $E/K$.  We can thus define a
character $\psi_{K'_{\tau'}}$ by composition.
%, i.e., composing $\psi_K$ with
%the map $\gal( K'_{\tau'}E/K'_{\tau'}) \to \gal(E/K)$.
Now, given degree one primes $\pK \in \unramified_{1}(K)$ and
$\pKp \in \unramified_{1}(K'_{\tau'})$, with $\pKp|\pK$, we will need,
no matter which  $\pKp | \pK$ is chosen, that
$$
\psi_{K'_{\tau'}}( \pKp ) = \psi_{K}( \pK ).
$$
%i.e., that the choice of $\pKp|\pK$ does not matter.
This in turn is
immediate from  \cite[p. 198, property A3]{LangANT}: since $\pKp$ and
$\pK$ are degree one primes, we have $f_{\pKp/\pK} = 1$, and thus
% $$
% \operatorname{res}_{E}(\pKp , K'_{\tau'}E/K'_{\tau'}) =
% (\pK, E/K)^{f_{\pKp/\pK}} = (\pK, E/K)
% $$
$$
\artin{K'_{\tau'}E/K'_{\tau'}}{\pKp} \bigg|_{E} =
\artin{E/K}{\pK}^{f_{\pKp/\pK}} =
\artin{E/K}{\pK}
$$
which shows that the choice of $\pKp|\pK$ indeed does not matter.

% Since each prime $\pLp \in \unramified(L')$ lies above a unique prime
% $\pK \in \unramified(K)$ we obtain a map $\pLp \to \pK$ from
% $\unramified(L')$ to $\unramified(K)$. We claim that the restriction
% of this map, from the set
% $$
% \{ \pLp \in \unramified(L')_{n' (??)}: \artin{L'/K}{\pLp} 
% = \tau'\}_{\tau' \in C'}
% $$
% to the set
% $$
% \{ \pK \in \unramified_{1}(K) : \artin{L'/K}{\pK}  = C'\},
% $$
% is $m:1$. To see this, note that given two primes
% $\pLp_{1},\pLp_{2} \in \unramified(L')$ that lie above
% $\pK \in \unramified(K)$, we have
% $\artin{L'/K}{\pLp_{1}} \sim \artin{L'/K}{\pLp_{2}}$ (where $\sim$
% denotes ``is conjugate to'').  We next ``descend'' from primes in $L'$
% down to primes in $K'_{\tau'}$; this map is $1:1$ as each $\pLp$ lies above a
% unique prime $\pKp$ in $K'_{\tau'}$ (under the assumption that
% $\artin{L'/K}{\pLp} \in C'$.)

Now consider the tower of extensions $L'/K'_{\tau'}/K$.  As for their
field of constants, we have $L' \cap \overline{\fq} = \fqnp$,
$K'_{\tau'} \cap \overline{\fq} = \fq$, $K \cap \overline{\fq} = \fq$.
Let $\pLp \in \unramified(L')$,
$\pKp \in \unramified_{1}(K'_{\tau'})$, and $\pK \in C_{1}(L/K, C)$
denote primes such that $\pLp | \pKp | \pK$.  Then, as
$\deg(\pKp) = \deg( \pK) = 1$ (see
\cite[Lemma~6.4.2]{FJ} and the proof of
Proposition~6.4.8) we have $\OKp/\pKp \simeq \OK/\pK
= \fq$, and thus 
$f_{\pKp/\pK} = 1$.  Further, as
$\artin{L'/K'_{\tau'}}{\pKp} = \{\tau'\}$ and $\ord(\tau') = n'$, we
have $\OLp/\pLp \simeq \fqnp$ and hence $f_{\pLp/\pKp} = n'$.  Thus,
as $[L':K'_{\tau'}] = n'$, the latter implies that any $\pLp | \pKp$
is uniquely determined by $\pKp$, whereas the former (together with
$[K'_{\tau'}:K]=m$) implies that there are exactly $m$ primes
$\pKp_{1}, \ldots, \pKp_m$ lying above $\pK$.

We thus find that
\begin{equation}
  \label{eq:sum-upstairs}
%\sum_{\pK \in C_{1}(L/K, C)} \psi_{K}(\pK)
%=
\sum_{\pK \in C_{1}(L'/K, C')} \psi_{K}(\pK)
=
\sum_{\tau' \in C'}
\sum_{\pK \in \unramified_{1}(K)}  \psi_{K}(\pK)
\sum_{ \substack{ \pLp \in \unramified_{n'}(L'), \pLp | \pK,  \\
    \artin{L'/K}{\pLp} = \{\tau' \}}} \frac1m
\end{equation}
$$
=
\frac1m 
\sum_{\tau' \in C'}
\sum_{\substack{ \pKp \in \unramified_{1}(K'_{\tau'}) \\
    \artin{L'/K'_{\tau'}}{\pKp} = \{\tau'\}}}
\psi_{K'_{\tau'}}(\pKp)
$$
which, from the proof of
\cite[Lemma~6.4.6]{FJ}, equals
$$
\frac1m 
\sum_{\tau' \in C'}
\sum_{\substack{ \pKp \in \unramified_{1}(K'_{\tau'}) }}
 \psi_{K'_{\tau'}}(\pKp);
$$
note that the inner sum amounts to summing an additive character over
{\em all} points on the curve given by the field $K'_{\tau'}$, having
field of constants
$K'_{\tau'} \cap \overline{\fq} = \fq$.

As $|C| = O_{[L:K]}(1)$, it suffices to show that
\begin{equation}
  \label{eq:curve-sum}
\left|
\sum_{\substack{ \pKp \in \unramified_{1}(K'_{\tau'}) }}
\psi_{K'_{\tau'}}(\pKp)
\right|
\ll_{[L:K]} q^{1/2}.
\end{equation}

To see this, let
$L(s) = L(s,\psi_{K'_{\tau'}}, EK'_{\tau'}/K'_{\tau'})$ denote the
Artin $L$-function attached to the character $\psi_{K'_{\tau'}}$; as
$ EK'_{\tau'}/K'_{\tau'}$ is abelian it is in fact an $L$-series
attached to a Hecke character.

Taking logarithmic derivatives of $L(s)$, we find that the sum over
primes in \eqref{eq:curve-sum} agrees, apart from ramified primes,
with the sum over degree one terms in $-L'(s)/L(s)$.  As $E/K$
only ramifies (wildly) at infinity, there are $O_{[L:K]}(1)$ ramified
primes, and hence the Riemann hypothesis for curves (due to Weil)
gives that the sum is $\ll_{[L:K]} q^{1/2}$, \emph{provided} that we can show
that the degree of the $L$-function only depends on $[L:K]$ (here some
care is needed since $[EK'_{\tau'}:K'_{\tau'}]=p$.)  To see this, first
note that
$$
\zeta_{K'_{\tau'}}(s) \cdot
\prod_{\psi \in  \widehat{\gal(EK'_{\tau'}/K'_{\tau'})} \setminus 1}
L(s,\psi, EK'_{\tau'}/K'_{\tau'})
=
\zeta_{EK'_{\tau'}}(s).
$$
Next note that any two non-trivial characters  $\psi_{1},\psi_{2} \in
\widehat{\gal(EK'_{\tau'}/K'_{\tau'})}$ have the same conductor, thus
the degrees of the $L$-functions
$$
\{ L(s,\psi, EK'_{\tau'}/K'_{\tau'})
\}_{\psi \in \widehat{\gal(EK'_{\tau'}/K'_{\tau'})} \setminus 1}
$$
are the same (cf. \cite[Ch.~6]{LiBook}.)
As there are $p-1$ non-trivial characters $\psi$, it is
enough to show that the degree of $\zeta_{EK'_{\tau'}}(s)$ is
$\ll_{[L:K]} p$.  As this degree is linear in $g = g(EK'_{\tau'})$,
the genus of $EK'_{\tau'}$, it is in turn enough to show that
$g \ll_{[L:K]} p$.  This is immediate from Castelnuovo’s Inequality
(cf. \cite[Theorem~3.11.3]{Stichtenoth-book}) since $g(E) \ll 1$,
$[E:K] = p$, $g(K'_{\tau'}) \ll_{[K:L]} 1$, and $[K'_{\tau'}:K]
\ll_{[L:K]} 1$.

% \section{Criteria for Morse polynomials and maximal Galois group}

\section{Criteria for Morse polynomials and  Galois group $S_{d}$}
\label{sec:crit-morse-polyn}

Given an integer $m \in [1,d-1]$, define
$$ 
f_{t}(x) := f(x) + t x^{m}.
$$

As long as a few of the small degree coefficients of $f$ avoid sets of
$O_{d}(1)$ ``bad coefficients'', the polynomial $f_{t}(x) $, over
$\fp(t)$, will have maximal Galois group.  To see this we recall a very
useful criterium.
\begin{prop}[Geyer \cite{JarRaz00}]
Assume that $\varphi = f/g \in \fp(x)$  is a Morse function of degree $d=
\deg(f) > \deg(g) \geq 0$.  Then the Galois group of the covering $\varphi :
{\mathbb P}^{1} \to {\mathbb P}^{1}$ of degree $d$ is the full
symmetric group, i.e., $\gal(f(x) - t g(x) / \fp(t)) \simeq S_{d}$.
\end{prop}

Fixing $g$, Geyer in fact shows that the set of polynomials
$f(x) = x^{d} + a_{d-1}x^{d-1}+ \cdots + a_{1}x + a_{0} $ for which
$f/g$ is Morse is a Zariski open dense subset of the affine $d$-space
with coordinates $a_{0}, \ldots, a_{d-1}$.
More precisely, for $\deg(g) > 0$, fixing $a_{d-1}, \ldots, a_{2}, a_{1}
\in \overline{\fp}$, write
$$
f(x) = x^{d} + a_{d-1} x^{d-1} + \cdots + a_{2} x^{2} + a_{1} x + u =
f_{\circ}(x) + u
$$
with $u$ trancendental over $\fp$.  Geyer then shows that
$f/g$ is Morse provided that
\begin{equation}
  \label{eq:1}
\gcd( f_{\circ}', g' ) = 1, \quad f_{\circ}'' \neq 0.
\end{equation}
In particular, if the conditions in (\ref{eq:1}) are satisfied, then
for all but $O_{d}(1)$ ``bad'' specializations of
$u=a_{0} \in \overline{\fp}$, the specialized ratio $f/g$ will be
Morse.

In particular, for $p$ large the second condition is automatic, and by
varying the linear coefficient of $f$ (again avoiding $O_{d}(1)$ ``bad''
values) we can ensure that the first condition $(f_{\circ}',g') = 1$ holds,
showing that $f(x) + t \cdot g(x)$ very often has full Galois group.
(Note that $f(x) + t \cdot g(x)$ and $f(x) - t \cdot g(x)$ have the
same Galois group over $\fp(t)$.)

We remark that the Morse criterion is certainly not needed for the
Galois group to be maximal.
E.g. (cf. \cite[Ch.~4.4]{serre-topics-in-galois-theory-book}) 
we have
$$
\operatorname{Gal}(x^{d} - x^{d-1} -t / \fp(t)) \simeq S_{d}.
$$

\subsection{Proof of Theorem~\ref{thm:primes-in-incomplete-interval}}
\label{sec:proof-theorem}

In \cite{kr-shortintervals} it was shown that in the family of
polynomials $f^{s}(x) = f(x) + 
sx$, $f^{s}$ is Morse for all but $O_{d}(1)$ values of $s \in \fp$ (or
even $s \in \overline{\fp}$.)  In particular, defining $B_{1}$ as the
set of $s$-values for which $f^{s}$ is not Morse, we have $|B_{1}| =
O_{d}(1)$ and the first part of
Theorem~\ref{thm:primes-in-incomplete-interval} follows from
Theorem~\ref{thm:incomplete-chebotarev} since $\gal(f(x)-t/\fp(t))
\simeq S_{d}$ for $s \in \fp \setminus B_{1}$.

Given Geyer's criterion, together with
Theorem~\ref{thm:incomplete-chebotarev}, the rest of the proof of
Theorem~\ref{thm:primes-in-incomplete-interval} is a simple matter
of checking the above conditions.

First, fix an integer $m \in [2,d-1]$ and take $g(x) = x^{m}$.  Then,
for $a_{1} \neq 0$, we have $(f_{\circ}', g') = 1$; as long as
$p > d+1$ we have $f_{\circ}'' \neq 0$, and hence, for all but
$O_{d}(1)$ choices of $a_{0} \in \overline{\fp}$ (and $a_{1}\neq 0$),
$f/g$ is Morse and
$$
\gal( f(x) + t x^{m} /  \fp(t)) \simeq S_{d}.
$$
Letting $B_{3}$ be the union, over $m \in [2,d-1]$, of these sets of
$O_{d}(1)$ exceptional $a_{0}$-values, the third part follows.

Similarly, for $m=1$ and $g(x) = x^{m} = x$ we may, possibly after
replacing $t$ by $t+1$ (changing the interval $I$ to an interval with the same cardinality, and symmetric difference of cardinality 2 with $I$), assume that $a_{1} \neq 0$.  Hence for all but
$O_{d}(1)$ choices of $a_{0}$, we find that
$\gal(f(x) + tx + a_{0} / \fp(t)) \simeq S_{d}$, and the second part
follows.

\subsection{Further examples}
\label{sec:further-examples}

We can also give examples of families of polynomials $f(x)$, with
fairly large number of free parameters (about $d/2$ of them),
such that the (geometric) Galois group of $f(x) + t x^{m}$ is not the
full symmetric group.  For instance, with $m=0$ and $d$ even, take
$f_{t}(x) = x^{d} + \sum_{i=0}^{d/2} a_{2i} x^{2i} + t$; the geometric
galois group is then a subgroup of a certain wreath product (here the
crucial point is that $f_{t}(x) = g_{t}(x^{2})$, i.e., the family is
decomposable.)  Similarly, for $m > 0$ and $m | d$, the family
$$
f_{t}(x) = x^{d} + \sum_{i=0}^{d/m} a_{i} x^{mi} + t x^{m} = g_{t}(x^{m})
$$
is decomposable.

%\bibliographystyle{abbrv} 
%\bibliography{bibshortint}

\end{document}